\documentclass[
hidelinks,onefignum,onetabnum,onealgnum]{siamart251216}

\usepackage[nocompress]{cite}
\usepackage{enumitem}
\usepackage{amssymb}
\usepackage{algpseudocode}
\usepackage{mathtools}
\usepackage{subcaption}
\usepackage{booktabs}
\usepackage{dsfont}

\newsiamthm{assumption}{Assumption}
\newsiamremark{remark}{Remark}
\newsiamremark{notation}{Notation}

\AddToHook{env/assumption/begin}{\crefalias{theorem}{assumption}}
\AddToHook{env/remark/begin}{\crefalias{theorem}{remark}}
\AddToHook{env/notation/begin}{\crefalias{theorem}{notation}}

\crefname{assumption}{Assumption}{Assumptions}
\Crefname{assumption}{Assumption}{Assumptions}
\crefname{remark}{Remark}{Remarks}
\Crefname{remark}{Remark}{Remarks}
\crefname{notation}{Notation}{Notations}
\Crefname{notation}{Notation}{Notations}

\newcommand{\scn}{{\scriptscriptstyle \mathcal{N}}}
\newcommand{\scDelta}{{\scriptscriptstyle \Delta}}
\newcommand{\N}{\mathbb{N}}
\newcommand{\R}{\mathbb{R}}
\newcommand{\C}{\mathbb{C}}
\newcommand{\restr}[2]{#1\!\restriction_{#2}}

\newcommand{\dY}[2]{
	\langle #1 , #2 \rangle _Y}
\newcommand{\dW}[2]{
	\langle #1 , #2 \rangle _W}

\newcommand{\dWs}[2]{
	\langle #1 , #2 \rangle _{W^*}}

\newcommand{\dWKs}[2]{
	\langle #1 , #2 \rangle _{W^*(K)}}
\newcommand{\dWK}[2]{
	\langle #1 , #2 \rangle _{W(K)}}
\newcommand{\duh}[2]{
	\langle #1 , #2 \rangle _h}
\newcommand{\duhm}[2]{
	\langle #1 , #2 \rangle _{h,\mu}}

\newcommand{\iK}[2]{
	(#1, #2)_K}
\newcommand{\iO}[2]{
	(#1, #2)_\Omega}

\newcommand{\iOO}[2]{
	(#1, #2)_{\Omega_0}}

\newcommand{\iY}[2]{
	(#1, #2)_Y}
\newcommand{\iWhs}[2]{
	(#1, #2)_{W_h^*}}

\newcommand{\iX}[2]{
	(#1, #2)_{X}}

\newcommand{\nK}[1]{
	\lVert #1 \rVert_K}
\newcommand{\nO}[1]{
	\lVert #1 \rVert_\Omega}

\newcommand{\nOO}[1]{
	\lVert #1 \rVert_{\Omega_0}}
\newcommand{\nWK}[1]{
	\lVert #1 \rVert_{W(K)}}
\newcommand{\nWhs}[1]{
	\lVert #1 \rVert_{W_h^*}}
\newcommand{\nWhsp}[1]{
	\lVert #1 \rVert_{(W_h^*)'}}

\newcommand{\nWO}[1]{
	\lVert #1 \rVert_{W}}
\newcommand{\nWOm}[1]{
	\lVert #1 \rVert_{W_\mu}}

\newcommand{\nWs}[1]{
	\lVert #1 \rVert_{W^*}}
\newcommand{\nQ}[1]{
	\lVert #1 \rVert_Q}
\newcommand{\nQm}[1]{
	\lVert #1 \rVert_{Q_\mu}}

\newcommand{\nWKs}[1]{
	\lVert #1 \rVert_{W^*(K)}}
\newcommand{\nUtQ}[1]{
	\lVert #1 \rVert_{H \times Q}}
\newcommand{\nX}[1]{
	\lVert #1 \rVert_X}

\newcommand{\nXm}[1]{
	\lVert #1 \rVert_{X_\mu}}
\newcommand{\nY}[1]{
	\lVert #1 \rVert_Y}
\newcommand{\nYp}[1]{
	\lVert #1 \rVert_{Y'}}
\newcommand{\nYK}[1]{
	\lVert #1 \rVert_{Y(K)}}
\newcommand{\nYm}[1]{
	\lVert #1 \rVert_{Y_\mu}}
\newcommand{\nYDp}[1]{
	\lVert #1 \rVert_{(Y_\scn^\scDelta)'}}

\newcommand{\infsupXY}[4]{
	\adjustlimits\inf_{0 \neq #1 \in #3} \sup_{0 \neq #2 \in #4} \frac{|b(#1,#2)|}{\nX{#1}\nY{#2}}}

\newcommand*\diff{\mathop{}\!\mathrm{d}}
\newcommand*\ii{\mathrm{i}}

\newcounter{assumption}

\newcounter{problem}
\renewcommand{\theproblem}{P\arabic{problem}}
\newcounter{schrproblem}

\newcommand{\vtick}{\mathord{\hspace{0.7pt}\raisebox{0.15ex}{\rule{0.47pt}{0.9ex}}}}

\newcommand{\tickplus}{\mathord{%
		\raisebox{0.15ex}{\hbox{%
				\rlap{\hbox to 0.90ex{\hfil\raisebox{0.45ex}{\rule{0.90ex}{0.45pt}}\hfil}}%
				\hbox to 0.90ex{\hfil\rule{0.46pt}{0.90ex}\hfil}%
		}}%
}}

\headers{Model order reduction for space-time Schrödinger}{T. Keith and M. Ohlberger}

\title{Model order reduction for space-time DPG approximations of the Schrödinger equation\thanks{Funded by the Deutsche Forschungsgemeinschaft (DFG, German Research Foundation) under Germany's Excellence Strategy EXC 2044/2 –390685587, Mathematics Münster: Dynamics–Geometry–Structure.}}

\author{Timo Keith\thanks{Institute for Analysis and Numerics, Faculty of Mathematics and Computer Science, University of Münster, Germany 
		(\email{timo.keith@uni-muenster.de}, \email{mario.ohlberger@uni-muenster.de}).}\and Mario Ohlberger\footnotemark[2]}

\ifpdf
\hypersetup{
	pdftitle={Model order reduction for space-time DPG approximations of the Schrödinger equation },
	pdfauthor={T. Keith, M. Ohlberger}
}
\fi

\begin{document}

\maketitle

\begin{abstract}
	This work extends the ultraweak space-time discontinuous Petrov--Galerkin (DPG) method for the linear time-dependent Schrödinger equation, introduced in [Demkowicz, Gopalakrishnan et al., SINUM 2017], from the free equation to bounded, time-independent, real-valued potentials. Continuous well-posedness is established in the graph norms with a potential-independent inf-sup constant. 
	The potential is parametrized and the discretization is combined with model order reduction. To the best of the authors' knowledge, this is the first explicit coupling of the DPG methodology with reduced basis methods. Uniform well-posedness with respect to parameter-independent reference norms is proven and yields exponential decay of the Kolmogorov $N$-widths. Moreover, the distinctive structural features of this combination are discussed, in particular the parameter dependence of the optimal test spaces and norms, and the built-in residual-based error control of the DPG framework, which drives a greedy reduced basis construction. 
	Numerical experiments in one spatial dimension for both parameter-independent and parametrized potentials support the theoretical results and demonstrate the accuracy and efficiency of the resulting reduced-order models.
\end{abstract}








\section{Introduction}
\label{sec:intro}
The time-dependent Schr\"odinger equation (TDSE) is one of the fundamental governing equations in quantum mechanics and describes the temporal evolution of a quantum-mechanical system \cite{quantMath}. For a single particle subject to a bounded external potential $V$ and, more generally, an inhomogeneous source term $f$, the TDSE, complemented by homogeneous Dirichlet boundary and initial conditions, reads as
\begingroup
\renewcommand{\theequation}{S0}
\begin{subequations}\label{S0-all}
	\renewcommand{\theequation}{\theparentequation.\arabic{equation}}
	\begin{alignat}{3}
		\label{Schrödingerclass} \ii \partial_tu - \Delta_x u + Vu &= f, &\qquad &x \in 	\Omega_0, &\ &t \in (0,T), \\
		\label{boundarycond} u(x,t) &= 0, &\qquad &x \in \partial \Omega_0, &\ &t 	\in (0,T),\\
		\label{initialcond} u(x,0) &= 0, &\qquad &x \in \Omega_0.
	\end{alignat}
\end{subequations}
\endgroup
Here, the set $\Omega_0 \subseteq \R^n$ denotes the spatial domain and $T > 0$ a finite time horizon.

The vast majority of established numerical methods for the TDSE follows a \emph{method-of-lines} approach, discretizing the spatial domain first and advancing the resulting system of ordinary differential equations in time by a suitable integrator. Full space-time discretizations have gained growing interest in recent years, particularly due to their ability to capture low-regularity solutions and their combination with \emph{model order reduction} techniques. For the Schr\"odinger equation, they are relatively scarce and include the Galerkin time-stepping methods of Karakashian and Makridakis for the cubic nonlinear equation \cite{spaceTimeNonlinarSEContin, spaceTimeNonlinarSEDiscontin} and the Trefftz and quasi-Trefftz discontinuous Galerkin methods of G\'omez and Moiola \cite{SpacetimeTrefftzDG, spacetimeDGVariablePotential}. Two global least-squares approaches are the isogeometric method of Bressan et al.\ \cite{spacetimeLSIMPORTANT} for the free equation and the framework of Dupuy et al.\ \cite{SpacetimeVFManyBody} for many-body systems. Very recently, Chakraborty et al.\ \cite{nonlinearSchrDPGSpaceTime} proposed a space-time \textit{discontinuous Petrov--Galerkin} (DPG) method for a modified model of laser pulse propagation.

Our interest lies in \textit{ultraweak} formulations, which transfer all derivatives to the test functions. Hain and Urban \cite{ultrWeakSchrodinger} developed an ultraweak framework that fixes the test space first and constructs the trial space accordingly, which can be interpreted as \emph{dual} to DPG. The direct foundation of this work is the space-time DPG formulation of Demkowicz, Gopalakrishnan et al.\ \cite{DPGSchr}, which works directly with the second-order Schr\"odinger operator. The authors show that the naive first-order reformulation with the flux $\tau \coloneqq \nabla_x u$ is \emph{not well-posed} in $L^2$. The second-order formulation requires no trace theory, but its well-posedness has to be established directly in the non-standard graph spaces of the space-time operator, which was found to be difficult for various space-time problems \cite[Section 11]{DPGActaNumerica}. So far, the concrete Schr\"odinger realization of \cite{DPGSchr} is restricted to the \emph{free} equation without potential. Extending it to a class of potentials broadens its applicability and, in particular, enables the parametrization of the potential, which opens the door to model order reduction.

In many applications, the problem depends on a parameter vector $\mu$ and has to be solved for many parameter values, for instance in parameter estimation, uncertainty quantification, or optimal control, where a full-order space-time solve for every parameter is infeasible. The \emph{reduced basis method} (RBM) \cite{RBMUrban} replaces the full-order model (FOM) by a reduced-order model (ROM) built from a few \textit{snapshots}, which a greedy algorithm selects by means of an a~posteriori error estimator. Traditional RBMs for time-dependent problems compute snapshot trajectories by time-marching schemes and build the basis, e.g., by the POD-Greedy algorithm \cite{PODGreedy}, but the resulting a~posteriori error bounds grow exponentially in time \cite{RBMTimeDependent}. Space-time variational formulations bypass this limitation, since treating time as an additional variable provides sharp, global-in-time error bounds. For ultraweak and optimally stable Petrov--Galerkin formulations, reduced basis schemes have been developed for transport problems \cite{doGre, brunkenTranspEq, MORUltraweakTransport} and, based on the \emph{dual} formulations of \cite{ultrWeakWave, ultrWeakSchrodinger}, for the wave and Schr\"odinger equations \cite[Section 1.5]{RBMUrban}. To the best of our knowledge, however, the DPG framework has not been explicitly combined with reduced basis techniques for any partial differential equation, although the connection is natural, since DPG methods already incorporate optimal test spaces and a built-in residual-based error estimator. Beyond these, the closest related works are reduced basis methods for minimum-residual mixed \cite{minResMixedRB} and least-squares \cite{LSFEMRB} finite element formulations of \textit{elliptic} problems, to which the DPG method is related as a minimum-residual method in a dual norm \cite[Theorem 3.2]{DPGActaNumerica}.

The contributions of this paper are the following. In \cref{sec:schroedinger}, we extend the ultraweak space-time DPG formulation of \cite{DPGSchr} to bounded, time-independent, real-valued potentials $V \in L^\infty(\Omega_0;\R)$ and prove continuous well-posedness in the graph norms with stability constants that are uniform in the potential, which is precisely what the parametric setting requires. For a parametrized potential, we establish in \cref{sec:mor} uniform well-posedness with respect to parameter-independent reference norms and exponential decay of the Kolmogorov $N$-widths of the solution manifold, a question left open for the dual formulations in \cite{RBMUrban}. We then develop the general combination of DPG methods with the reduced basis method. The optimal test spaces yield a feasible supremizer-based reduced test space, and the built-in residual provides an error estimator with respect to the infinite-dimensional solution, which drives a greedy basis construction. Finally, the scheme is realized in a hybrid implementation, static condensation for the full-order snapshots and normal equations for the reduction, and validated in one spatial dimension in \cref{sec:numerics}. Appendix \ref{sec:dpg} collects the background on the DPG method.

\section{The Schrödinger equation with a spatially varying potential}
\label{sec:schroedinger}
This section extends the continuous well-posedness theory of \cite{DPGSchr} to bounded, real-valued, time-independent potentials, introducing the constructions therein, particularly of its appendix, directly for the Schrödinger case. Well-posedness of the ultraweak formulation follows from that of an auxiliary \textit{strong} formulation, in whose proof the potential enters only through the spectral decomposition of the spatial operator. The key observation is an a~priori bound for this strong formulation whose constant depends on the time horizon $T$ only, so that the inf-sup constant of the ultraweak formulation is bounded from below uniformly in $V$, which is what  its usage in \cref{sec:mor} requires.

In the following, we consider \cref{S0-all} for some bounded, open Lipschitz domain $\Omega_0 \subseteq \R^n$ and potential $V \in L^\infty(\Omega_0;\R)$, and define the space-time domain $\Omega \coloneqq \Omega_0 \times (0,T) \subseteq \R^d$, $d \coloneqq n + 1$. Let $\Omega_h$ be some partition of $\Omega$ into a mesh of finitely many, nonempty open elements $K \in \Omega_h$ satisfying $\cup_{K \in \Omega_h} \bar{K} = \bar{\Omega}$, with $h \coloneqq \max_{K \in \Omega_h} \text{diam}(K)$.
\subsection{Functional setting}
\label{sec:setting}
We set the boundary parts
\begin{align}
	\label{Gamma}
	\Gamma \coloneqq \partial \Omega_0 \times [0,T] \cup \Omega_0 \times \{0\}, \qquad \Gamma^* \coloneqq \partial \Omega_0 \times [0,T] \cup \Omega_0 \times \{T\}.
\end{align}
For any nonempty, open subset \(K\subseteq\Omega\), $\mathcal{D}(K)\coloneqq C_c^\infty(K)$ and $\mathcal{D}'(K)$ denote the space of smooth compactly supported functions and its dual, and $\mathcal{D}(\bar{\Omega}) \coloneqq \{\restr{\varphi}{\Omega} \, \mid \, \varphi \in \mathcal{D}(\R^d)\}$ the restrictions of smooth compactly supported functions on $\R^d$ to $\Omega$, analogously $\mathcal{D}(\bar{\Omega}_0)$. We denote the $L^2$-inner product and $L^2$-norm on $K$ by $\iK{\cdot}{\cdot}$ and $\nK{\cdot}$, respectively. In the remainder of this work, $L^2$-, $C^k$-, and $\mathcal{D}$-spaces always refer to \textit{complex-valued} functions, inner products are antilinear in their second argument, and primes denote conjugate-duals.
Just as in \cite[Problem 2.3]{DPGSchr}, we first aim at a strong variational formulation of \cref{S0-all}. 
\begin{equation}
	\label{strongform}
	\stepcounter{problem} \tag{\theproblem}
	\begin{aligned}
		\text{For some given } f \in L^2(\Omega), \text{ find } u \in U \text{ satisfying}\\ 
		\iO{Au}{v} = \iO{f}{v} \ \text{for all } v \in L^2(\Omega). \qquad
	\end{aligned}
\end{equation}
Here,
\begin{alignat}{3}
	\label{ASchr}
	A&: &\, &L^2(\Omega) \to \mathcal{D}'(\Omega), &\quad A \coloneqq& \, \ii \partial_t - \Delta_x + V, \\
	\label{AstarSchr}
	A^*&: &\, &L^2(\Omega) \to \mathcal{D}'(\Omega), &\quad A^* =& \, \ii \partial_t - \Delta_x + V,
\end{alignat}
with the formal adjoint $A^*$, which coincides with $A$ since $V$ is real-valued. The domain $U \subseteq L^2(\Omega)$ is defined via the (local) graph spaces and boundary operators from the appendix of \cite{DPGSchr}. 
\begin{definition}
	\label{def:graphspacesboundaryoperators}
	We define the \textit{graph spaces}
	\begin{alignat}{3}
		\label{WK}
		W(K) &\coloneqq \{ u \in L^2(K) \ | \ Au \in L^2(K)\}, &\quad &\nWK{\cdot} &\coloneqq &\, (\nK{\cdot}^2 + \nK{A \cdot}^2)^{1/2}, \\
		\label{WKs}
		W^*(K) &\coloneqq \{ v \in L^2(K) \ | \ A^* v \in L^2(K)\}, &\quad &\nWKs{\cdot} &\coloneqq &\, (\nK{\cdot}^2 + \nK{A^* \cdot}^2)^{1/2},
	\end{alignat}
	which are Hilbert spaces by \cite[Lemma A.1]{DPGSchr}, set $W \coloneqq W(\Omega)$ and $W^* \coloneqq W^*(\Omega)$, and define the bounded (local) \textit{boundary operators} for $u \in W(K)$, $v \in W^*(K)$ by
	\begin{alignat}{2}
		\label{DK}
		D_K: W(K) \to W^*(K)'&, &\quad \dWKs{D_Ku}{v} &\coloneqq \iK{Au}{v} - \iK{u}{A^* v},\\
		\label{DKs}
		D_K^*: W^*(K) \to W(K)'&, &\quad \dWK{D_K^*v}{u} &\coloneqq \iK{A^*v}{u} - \iK{v}{Au},
	\end{alignat}
	with global versions $D \coloneqq D_\Omega$ and $D^* \coloneqq D^*_\Omega$. In the present case, $W(K) = W^*(K)$ and $D_K = D_K^*$ coincide, as $A = A^*$ do. Adapting the notation of \cite{DPGSchr} ($V$ denoting the potential here), we define
	\begin{align}
		\label{mathcalU}
		\mathcal{U} &\coloneqq \{\varphi \in \mathcal{D}(\bar{\Omega}) \ | \ \restr{\varphi}{\Gamma} = 0\}, \qquad
		\mathcal{U}^* \coloneqq \{\varphi \in \mathcal{D}(\bar{\Omega}) \ | \ \restr{\varphi}{\Gamma^*} = 0\}, \\
		\label{domASchr}
		U &\coloneqq \{u \in W \ | \ \dW{D^*v}{u} = 0 \ \forall\, v \in \mathcal{U}^*\} = {^\perp D^*}(\mathcal{U}^*), 
	\end{align} 
	the latter being a closed subspace of $W$ equipped with the $\nWO{\cdot}$-graph norm, and the domain of the adjoint
	\begin{align*}
		U^* \coloneqq \{s \in L^2(\Omega) \ | \ \exists \, l \in L^2(\Omega) : \iO{Au}{s} = \iO{u}{l} \ \forall \,u \in U\} = {^\perp D}(U),
	\end{align*}
	equipped with the $\nWs{\cdot}$-graph norm, where the last identity is \cite[Lemma A.2]{DPGSchr}.
\end{definition}
\begin{proposition}
	\label{prop:potentialIndepGraphSpaces}
	All constructions from \cref{def:graphspacesboundaryoperators} are independent of the potential. More precisely, denoting by $\tilde{A} \coloneqq \ii\partial_t - \Delta_x$ the potential-free operator and by $\tilde{W}(K)$, $\tilde{D}_K$, $\tilde{U}$, $\tilde{U}^*$ the corresponding constructions, we have $W(K) = \tilde{W}(K)$ as sets with
	\begin{align}
		\label{eq:normequiv}
		(2 + \|V\|_\infty^2)^{-1/2} \, \nWK{u} \leq \|u\|_{\tilde{W}(K)} \leq (2 + \|V\|_\infty^2)^{1/2} \, \nWK{u} 
	\end{align}
	for all $u \in W(K)$, $D_K = \tilde{D}_K$ as operators, and $U = \tilde{U}$, $U^* = \tilde{U}^*$ as sets with equivalent norms.
\end{proposition}
\begin{proof}
	The graph spaces coincide as sets since the potential is bounded, and \eqref{eq:normequiv} follows from $\nK{Au - Vu}^2 \leq (1 + \|V\|_\infty^2) (\nK{Au}^2 + \nK{u}^2)$ by the Cauchy--Schwarz inequality in $\R^2$, with the roles of $A$ and $\tilde{A}$ interchanged for the reverse direction. The boundary operators are independent of the potential since $\iK{Vu}{w} = \iK{u}{Vw}$ for any $u, w \in W(K)$ by its real-valuedness, which also yields the claims for $U$ and $U^*$.
\end{proof}
As in \cite[Assumption 1]{DPGSchr}, well-posedness will rest on the following density assumption. 
\begin{assumption}
	\label{ass:densityU}
	The spaces $\mathcal{U} \subseteq U$ and $\mathcal{U}^* \subseteq U^*$ are dense with respect to the respective graph norms. 
\end{assumption}
By \cref{prop:potentialIndepGraphSpaces}, \cref{ass:densityU} does not depend on the potential, i.e., it holds for \textit{every} $V \in L^\infty(\Omega_0;\R)$ once it is known for the free equation (cf.\ \cref{thmdensity1DSchr}). Under \cref{ass:densityU}, \cite[Lemma 2.2 and (2.3), (2.4)]{DPGSchr} moreover yield
\begin{align}
	\label{eq:annihilators}
	U = {^\perp D^*}(U^*) \qquad \text{and} \qquad U^* = {^\perp D}(\mathcal{U}).
\end{align}

Next, let $W_h \coloneqq \prod_{K \in \Omega_h} W(K) = \prod_{K \in \Omega_h} W^*(K) \eqqcolon W_h^*$ denote the broken graph spaces with product norm $\nWhs{v}^2 \coloneqq \sum_{K \in \Omega_h} \nWKs{v}^2$, and let $D_h = D_h^*$ and $A_h = A_h^*$ be the above operators applied piecewise on each mesh element, i.e., $\duh{D_h w}{v} \coloneqq \iO{A_h w}{v} - \iO{w}{A_h^* v}$ for $w \in W_h$, $v \in W_h^*$, where $\duh{\cdot}{\cdot}$ denotes the duality pairing in $(W_h^*)'$. Note that $W \subseteq W_h$ continuously and $\duh{D_h w}{v} = \dWs{Dw}{v}$ for $w \in W$, $v \in W^*$ \cite[Lemma A.3]{DPGSchr}.
Testing \eqref{strongform} by $v \in W_h^*$ and integrating by parts element-wise yields $\duh{D_h u}{v} + \iO{u}{A_h^*v} = \iO{f}{v}$ for all $v \in W_h^*$, which motivates the following mesh-dependent ultraweak formulation.
\begin{definition}
	\label{def:Qb}
	We abbreviate $H \coloneqq L^2(\Omega)$, define the interface space 
	\begin{align}
		\label{Q}
		Q \coloneqq D_h(U) \subseteq (W_h^*)', \qquad \nQ{q} \coloneqq \inf_{w \in {\restr{D_h}{U}}^{-1} (\{q\})} \nWO{w},
	\end{align}
	which is complete with respect to this quotient norm \cite[Appendix A]{DPGSchr}, introduce the additional interface unknown $q \coloneqq D_h u \in Q$, and define the sesquilinear form
	\begin{align}
		\label{b}
		b: (H \times Q) \times W_h^* \to \C\ , \ b((u,q),v) \coloneqq \iO{u}{A_h^*v} + \duh{q}{v},
	\end{align}
	with $\nUtQ{(u,q)}^2 \coloneqq \nO{u}^2 + \nQ{q}^2$ as usual. For general $f \in (W_h^*)'$, the ultraweak formulation of \cref{strongform} reads as follows.
	\begin{equation}
		\begin{aligned}
			\text{Given some } f \in (W_h^*)' \text{, find } u \in H, \ q \in Q \ \text{ s.t.}\\
			b((u,q),v) = f(v) \ \ \forall \,v \in W_h^*. \qquad \qquad
		\end{aligned}
		\stepcounter{problem} \tag{\theproblem} \label{weakform}
	\end{equation}
\end{definition}
\subsection{Spectral decomposition of the spatial operator}
\label{sec:spectral}
The surjectivity argument of \cite[Theorem 2.4]{DPGSchr} separates variables with respect to a Dirichlet eigenbasis of $-\Delta_x$. With a potential, the corresponding ingredient is an eigenbasis of $-\Delta_x + V$. Its existence is classical for real-valued functions and, in the presence of a bounded potential, with the eigenvalue equation in variational form, see \cite[Theorem 9.31 and Remark 30]{Brezis}. We need the complex setting, the eigenvalue equation in $L^2(\Omega_0)$, and an integration by parts formula against smooth functions vanishing on $\partial \Omega_0$. We abbreviate $\|u\|_{1}^2 \coloneqq \nOO{u}^2 + \nOO{\nabla u}^2$ on $H_0^1(\Omega_0)$.
\begin{lemma}[Spectral decomposition of the spatial operator]
	\label{lem:spectral}
	There exist an orthonormal basis $\{e_k\}_{k = 1}^\infty$ of $L^2(\Omega_0)$ and real numbers $\omega_k$ such that, for every $k \in \N$, $e_k \in H_0^1(\Omega_0)$, its weak Laplacian $\Delta e_k \in L^2(\Omega_0)$ exists, and
	\begin{align}
		\label{eq:eigen-strong}
		-\Delta e_k + V e_k &= \omega_k e_k \qquad \text{a.e. in } \Omega_0,\\
		\label{eq:eigen-variational}
		\iOO{(-\Delta + V) e_k}{\psi} &= \iOO{e_k}{(-\Delta + V) \psi} \qquad \forall\, \psi \in \mathcal{D}(\bar{\Omega}_0) \text{ with } \restr{\psi}{\partial \Omega_0} = 0.
	\end{align}
\end{lemma}
\begin{proof}
	By \cite[Remark 30]{Brezis} applied with $a_{ij} = \delta_{ij}$ and $a_0 = V$, there exist an orthonormal basis $\{e_k\}_{k=1}^\infty$ of $L^2(\Omega_0;\R)$ and real numbers $\omega_k$ such that $e_k \in H_0^1(\Omega_0;\R)$ and
	\begin{align}
		\label{eq:eigen-weak}
		\iOO{\nabla e_k}{\nabla v} + \iOO{V e_k}{v} = \omega_k \iOO{e_k}{v} \qquad \forall\, v \in H_0^1(\Omega_0;\R).
	\end{align}
	Decomposing $f \in L^2(\Omega_0)$ into real and imaginary parts shows that the $e_k$ also form an orthonormal basis of the complex space $L^2(\Omega_0)$. Likewise, \eqref{eq:eigen-weak} holds for all $v \in H_0^1(\Omega_0)$ by sesquilinearity, since $\operatorname{Re} v, \operatorname{Im} v \in H_0^1(\Omega_0;\R)$.

	To prove \eqref{eq:eigen-strong}, let $\varphi \in \mathcal{D}(\Omega_0)$. The definition of the weak gradient of $e_k$, tested with $\partial_j \bar{\varphi}$, gives $\iOO{\nabla e_k}{\nabla \varphi} = - \int_{\Omega_0} e_k \, \Delta \bar{\varphi}$, so that \eqref{eq:eigen-weak} with $v = \varphi$ states precisely that $-\Delta e_k = \omega_k e_k - V e_k$ in $\mathcal{D}'(\Omega_0)$, as $\varphi \mapsto \bar{\varphi}$ maps $\mathcal{D}(\Omega_0)$ onto itself. Since the right-hand side lies in $L^2(\Omega_0)$, so does the weak Laplacian $\Delta e_k$, and \eqref{eq:eigen-strong} holds a.e.
	
	To prove \eqref{eq:eigen-variational}, let $\psi \in \mathcal{D}(\bar{\Omega}_0)$ vanish on $\partial \Omega_0$. As $\Omega_0$ is a bounded Lipschitz domain, $H_0^1(\Omega_0)$ is the kernel of the trace operator \cite[Theorem A8.6 and Lemma A8.10]{altLFA}, so that $\psi$ is an admissible test function in \eqref{eq:eigen-weak}. Furthermore, the integration by parts formula
	\begin{align}
		\label{eq:green}
		\iOO{\nabla u}{\nabla \psi} = \iOO{u}{-\Delta \psi} \qquad \forall\, u \in H_0^1(\Omega_0)
	\end{align}
	holds without boundary terms, because the boundary condition is carried by $u$. Indeed, for $u \in \mathcal{D}(\Omega_0)$ it is the classical formula, both sides are continuous in $u$ with respect to $\|\cdot\|_1$ since $\nabla \psi, \Delta \psi \in L^2(\Omega_0)$, and $\mathcal{D}(\Omega_0)$ is dense in $H_0^1(\Omega_0)$ by definition. Applying \eqref{eq:eigen-strong}, then \eqref{eq:eigen-weak} with $v = \psi$, and finally \eqref{eq:green} with $u = e_k$ together with the real-valuedness of $V$, we obtain
	\begin{align*}
		\iOO{(-\Delta + V) e_k}{\psi}
		&= \omega_k \iOO{e_k}{\psi} \\
		&= \iOO{\nabla e_k}{\nabla \psi} + \iOO{V e_k}{\psi}
		= \iOO{e_k}{(-\Delta + V) \psi}.
	\end{align*}
\end{proof}
\subsection{Well-posedness}
\label{sec:wellposed}
We are now in a position to generalize the well-posedness results of \cite{DPGSchr}, additionally providing an a~priori bound with a constant that does not depend on the potential.
\begin{theorem}[Strong formulation]
	\label{thmstrongSchrwellposed}
	Suppose that \cref{ass:densityU} holds. Then, $A: U \to L^2(\Omega)$ is a continuous bijection, so \cref{strongform} is well-posed, and
	\begin{align}
		\label{eq:apriori}
		\nO{u} \leq \frac{T}{\sqrt{2}} \nO{Au} \qquad \forall\, u \in U,
	\end{align}
	so that $\|A^{-1}\|_{L^2(\Omega) \to U} \leq C_T \coloneqq (1 + T^2/2)^{1/2}$, independently of $V$.
\end{theorem}
\begin{proof}
	The proof of \cite[Theorem 2.4]{DPGSchr} carries over verbatim with the eigenpairs $\{(e_k, \omega_k)\}_k$ of \cref{lem:spectral} in place of those of $-\Delta_x$ in the definitions \cite[(1.3)]{DPGSchr} of $f_k$, $u_k$, $F_M$, and $U_M$. Indeed, the eigenpairs enter only through the orthonormal basis property, through the identity $A U_M = F_M$, which follows from $\ii u_k' + \omega_k u_k = f_k$ and \cref{eq:eigen-strong}, through the integration by parts in space showing $U_M \in U$, which is precisely \cref{eq:eigen-variational} applied to $\varphi(\cdot,t)$ for $\varphi \in \mathcal{U}^*$, and through $|e^{\ii \omega_k (t-s)}| = 1$ in the Cauchy estimates, which holds since $\omega_k$ is real, positivity of the eigenvalues never being used. All remaining arguments concern the spaces $U$, $U^*$, $\mathcal{U}$, $\mathcal{U}^*$ and the boundary operator only, which are independent of the potential by \cref{prop:potentialIndepGraphSpaces} and satisfy \cref{eq:annihilators}. For \cref{eq:apriori}, the Cauchy estimate of \cite[Theorem 2.4]{DPGSchr} for $U_N - U_M$, whose termwise proof remains valid for $M = 0$ with $U_0 = 0$, reads $\nO{U_N}^2 \leq \frac{T^2}{2} \sum_{k=1}^{N} \int_0^T |f_k(t)|^2 \diff t \leq \frac{T^2}{2} \nO{f}^2$, hence $\nO{u} \leq \frac{T}{\sqrt{2}} \nO{f}$ for the constructed solution $u = \lim_{N \to \infty} U_N$. By injectivity, every $u \in U$ is the constructed solution for $f \coloneqq Au$, and $\nWO{u}^2 = \nO{u}^2 + \nO{Au}^2 \leq (1 + T^2/2) \nO{Au}^2$ yields the bound on $\|A^{-1}\|$.
\end{proof}
In contrast to \cite[Theorem 2.7]{DPGSchr}, we track the stability constants explicitly.
\begin{theorem}[Ultraweak formulation]
	\label{thmwellposedweakSchr}
	Suppose that \cref{ass:densityU} holds. Then, Problem~\eqref{weakform} is well-posed in the sense of \eqref{N1}-\eqref{N3} from Appendix \ref{sec:dpg} with continuity constant $\gamma = \sqrt{2}$ and inf-sup constant
	\begin{align}
		\label{eq:beta}
		\beta \coloneqq \bigl((C+1)^2 + C^2\bigr)^{-1/2} > 0, \qquad C \coloneqq \|A^{-1}\|_{L^2(\Omega) \to U} \leq C_T.
	\end{align}
	In other words, for any $f \in (W_h^*)'$, \eqref{weakform} admits a unique solution $(u,q) \in H \times Q$ with $\nUtQ{(u,q)} \leq \beta^{-1} \nWhsp{f}$. Moreover, if $f = \iO{g}{\cdot}$ for some $g \in L^2(\Omega)$, then $u \in U$ solves the strong formulation \eqref{strongform} and $D_h u = q$.
\end{theorem}
\begin{proof}
	Well-posedness and the regularity statement follow from \cite[Theorem A.5]{DPGSchr}, whose hypotheses we verify. Conditions (A-a) to (A-c) of \cite[Appendix A]{DPGSchr} hold since $A$ is a second-order operator with bounded coefficients and $\mathcal{D}(\Omega) \subseteq U$. Condition \cite[(A.12)]{DPGSchr}, i.e., $U = {^\perp D^*}(U^*)$, is \eqref{eq:annihilators}, condition \cite[(A.13)]{DPGSchr} is \cref{thmstrongSchrwellposed}. It remains to make the constants explicit.
	For \eqref{N1}, we calculate 
	\begin{align*}
		|b((u,q),v)| \leq \nO{u} \nO{A_h^*v} + \nWhsp{q}\nWhs{v} \leq \sqrt{2} \, \nUtQ{(u,q)} \nWhs{v},
	\end{align*}
	where $\nO{A_h^*v} \leq \nWhs{v}$ by definition of the broken graph norm, and $\nWhsp{q} \leq \nQ{q}$ since $|\duh{D_h w}{v}| \leq \nO{Aw}\nO{v} + \nO{w}\nO{A_h^* v} \leq \nWO{w} \nWhs{v}$ for all $w \in U$ with $D_h w = q$ (in fact, $\nWhsp{q} = \nQ{q}$ by \cite[Lemma A.10]{DPGSchr}). 
	For \eqref{N2}, let $v \in W_h^* \subseteq L^2(\Omega)$. By \cref{thmstrongSchrwellposed}, there exists some $w \in U$ with $A_hw = Aw = v$ and $\nWO{w} \leq C \nO{v} \leq C \nWhs{v}$. 
	Setting $q \coloneqq D_hw \in Q$ and $u \coloneqq w + A_h^*v \in H$, we calculate
	\begin{align*}
		\nWhs{v}^2 = \iO{A_hw}{v} + \iO{A_h^*v}{A_h^*v} 
		= \iO{w + A_h^*v}{A_h^*v} + \duh{D_hw}{v} 
		= b((u,q),v),
	\end{align*}
	where $\nQ{q} \leq \nWO{w} \leq C \nWhs{v}$ by \eqref{Q} and $\nO{u} \leq \nWO{w} + \nWhs{v} \leq (C+1) \nWhs{v}$, i.e., $\nUtQ{(u,q)} \leq \beta^{-1} \nWhs{v}$. This is the inf-sup condition \cite[(A.17)]{DPGSchr} with constant $\beta^{-1}$ and, since $B: H \times Q \to (W_h^*)'$ is an isomorphism, the same constant applies in \eqref{N2}.
\end{proof}
\begin{remark}[Uniformity of the constants]
	\label{rmk:constants}
	The inf-sup constant $\beta$ in \eqref{eq:beta} depends on the potential only through $C = \|A^{-1}\|$, and $C \mapsto ((C+1)^2 + C^2)^{-1/2}$ is decreasing. Hence, \cref{thmstrongSchrwellposed} yields
	\begin{align}
		\label{eq:betaT}
		\beta \geq \beta_T \coloneqq \bigl((C_T+1)^2 + C_T^2\bigr)^{-1/2}, \qquad C_T = (1 + T^2/2)^{1/2},
	\end{align}
	\textit{independently} of $V \in L^\infty(\Omega_0;\R)$, but with respect to the $V$-dependent graph norms. Relative to the potential-free norms, one pays the equivalence constants from \cref{prop:potentialIndepGraphSpaces}. Both facts will be the key to the uniform stability results in \cref{sec:param}.
\end{remark}
\begin{corollary}[One spatial dimension]
	\label{thmdensity1DSchr}
	Let $\Omega_0 = (0,L)$ for some $L > 0$. Then, \cref{ass:densityU} is satisfied, so that \cref{thmstrongSchrwellposed,thmwellposedweakSchr} hold unconditionally for every $V \in L^\infty((0,L);\R)$.
\end{corollary}
\begin{proof}
	\cite[Theorem 3.1]{DPGSchr} proves the density for the free equation by reflection, time translation, and smoothing. Since a potential changes neither the spaces $U$, $\mathcal{U}$, $U^*$, $\mathcal{U}^*$ as sets nor their topologies by \cref{prop:potentialIndepGraphSpaces}, the density transfers.
\end{proof}
\begin{remark}[Scope of the abstract framework]
	\label{rmkAssDensity}
	The reflection argument of \cite[Theorem 3.1]{DPGSchr} relies on the one-dimensional geometry. In the appendix of \cite{spacetimeLSIMPORTANT}, this argument is extended to unit hypercubes $\Omega_0 = (0,1)^n$ and certain \textit{isogeometric} domains, where \cref{thmstrongSchrwellposed,thmwellposedweakSchr} thus hold as well. For time-dependent potentials $V \in L^\infty(\Omega;\R)$, \cref{prop:potentialIndepGraphSpaces}, the potential-independence of \cref{ass:densityU}, and all subsequent discrete and reduced arguments remain valid verbatim. Only the bijectivity of $A: U \to L^2(\Omega)$ then needs a different proof, since the separation of variables in \cref{thmstrongSchrwellposed} relies on a time-independent spatial operator. The restriction to time-independent potentials thus stems from this step alone and not from the DPG machinery.
\end{remark}
\subsection{Discretization and a priori estimates in one spatial dimension}
\label{sec:discrete1D}
Let $H_\scn \subseteq H$ and $Q_\scn \subseteq Q$ be finite-dimensional subspaces, $X_\scn \coloneqq H_\scn \times Q_\scn$, and let \textnormal{\hypertarget{idealDPGSchr}{\text{(I-TDSE)}}} and \hypertarget{practicalDPGSchr}{\text{(P-TDSE)}} denote the ideal DPG method \eqref{DPG} and the practical DPG method \eqref{approxidealDPG} from Appendix \ref{sec:dpg} with an enriched test space $W_\scn^{*\scDelta} \subseteq W_h^*$ for \eqref{weakform}. Under \cref{ass:densityU}, \hyperlink{idealDPGSchr}{\text{(I-TDSE)}} is stable and quasi-optimal in the sense of \cref{thmideal} by \cref{thmwellposedweakSchr}, independently of the potential.

In one spatial dimension, $\Omega_0 = (0,L)$, we adopt the spaces of \cite[Section 4]{DPGSchr} on a uniform mesh $\Omega_h$ of space-time squares of side length $h$. For $p \in \N$, $H_\scn$ consists of the piecewise tensor-product polynomials of degree at most $p-1$. The interface space $Q_\scn$ is identified there with the product $Q_\scn^{\tickplus} \times Q_\scn^{\vtick}$ of continuous traces $q^{\tickplus}$ on all element edges, vanishing on $\Gamma$, and discontinuous fluxes $q^{\vtick}$, i.e., traces of the $x$-derivative, on the vertical edges, both of degree at most $p$. Together with the enriched test space $W_\scn^{*\scDelta}$ of piecewise tensor-product polynomials of degree at most $p + \Delta p$ for some $\Delta p \geq 1$, these spaces define \hyperlink{practicalDPGSchr}{\text{(P-TDSE)}} as used in \cref{sec:numerics}.

The a~priori estimate of \cite[Theorem 4.3]{DPGSchr} remains valid with a potential. It is restricted to one spatial dimension, as its proof approximates the interface variable through a conforming finite element subspace of $U$, whose functions have to be continuous in time and continuously differentiable in space, and \cite[Section 4]{DPGSchr} constructs such elements for $n = 1$ only.
\begin{theorem}[A~priori estimate]
	\label{thmerror1DSchrDPG}
	Let $p \geq 3$, $f = \iO{g}{\cdot}$ for some $g \in L^2(\Omega)$, and assume that the unique solution $(u,q) \in H \times Q$ of \eqref{weakform} satisfies $u \in H^{p+1}(\Omega)$. Then, for any $2 \leq r \leq p-1$, there exists a constant $C > 0$ independent of $h$ such that the unique solution $(u_\scn,q_\scn) \in H_\scn \times Q_\scn$ of \textnormal{\hyperlink{idealDPGSchr}{\text{(I-TDSE)}}} satisfies 
	\begin{align*}
		\nO{u - u_\scn} + \nQ{q-q_\scn} \leq C h^r |u|_{H^{r+2}(\Omega)}.
	\end{align*}
\end{theorem}
\begin{proof}
	As the potential-dependent graph norm is equivalent to the potential-independent one by \cref{prop:potentialIndepGraphSpaces}, the interpolation error estimates from \cite[Lemma 4.2]{DPGSchr} remain valid up to an additional $h$-independent constant depending on $\|V\|_\infty$, and the claim follows just as in \cite[Theorem 4.3]{DPGSchr}.
\end{proof}
\Cref{thmerror1DSchrDPG} concerns the ideal method, whereas all computations use \hyperlink{practicalDPGSchr}{\text{(P-TDSE)}} in the mixed form \eqref{mixedFormDiscrete} with $Y_\scn^\scDelta = W_\scn^{*\scDelta}$. The following two remarks concern its stability and its algebraic solution.
\begin{remark}[Discrete stability]
	\label{rmkDiscreteStab}
	For \hyperlink{practicalDPGSchr}{\text{(P-TDSE)}}, the Fortin condition \eqref{Fortin} has to be verified with respect to the second-order graph norm $\nWhs{\cdot}$. To the best of our knowledge, no such Fortin operator is available even for the free equation, the question not being addressed in \cite{DPGSchr, DPGActaNumerica}, and the known space-time constructions \cite{spacetimeDPGHeat, spacetimeAccWaves, spacetimeLinearWaveWieners} concern first-order settings. Hence, we rely on numerical evidence for the choice of $\Delta p$, as is common in the DPG literature, see \cref{sec:fom}.
\end{remark}
\begin{remark}[Algebraic realization]
	\label{rmkStaticCond}
	Splitting $b((u,q),v) = \iO{u}{A_h^* v} + \duh{q}{v}$, denoting by $\boldsymbol{\mathcal{G}}$ the Gramian of $W_\scn^{*\scDelta}$ with respect to $\iWhs{\cdot}{\cdot}$, by $\boldsymbol{\mathcal{B}}_H$ and $\hat{\boldsymbol{\mathcal{B}}}$ the matrices of the two parts of $b$, and by $\boldsymbol{e}$, $\boldsymbol{u}$, $\boldsymbol{q}$ the coefficient vectors of the solution components of \eqref{mixedFormDiscrete}, the mixed system takes the block form
	\begin{align}
		\label{eq:blocksystem}
		\begin{pmatrix}
			\boldsymbol{\mathcal{G}} & \boldsymbol{\mathcal{B}}_H & \hat{\boldsymbol{\mathcal{B}}} \\
			\boldsymbol{\mathcal{B}}_H^{H} & 0 & 0 \\
			\hat{\boldsymbol{\mathcal{B}}}^{H} & 0 & 0
		\end{pmatrix}
		\begin{pmatrix}
			\boldsymbol{e} \\ \boldsymbol{u} \\ \boldsymbol{q}
		\end{pmatrix}
		=
		\begin{pmatrix}
			\boldsymbol{\mathcal{F}} \\ 0 \\ 0
		\end{pmatrix}.
	\end{align}
	Since the volume component $H_\scn$ is broken as well, a feature specific to the ultraweak formulation, both $\boldsymbol{\mathcal{G}}$ and $\boldsymbol{\mathcal{B}}_H$ act on \textit{local} degrees of freedom, and \textit{static condensation} \cite{NGSolve} eliminates $(\boldsymbol{e}, \boldsymbol{u})$ by element-wise solves, leaving a Schur complement system for the interface unknowns $\boldsymbol{q}$ of size $\operatorname{dim}(Q_\scn)$ only. 
\end{remark}
\section{Model order reduction}
\label{sec:mor}
We now apply the reduced basis method \cite{RBMUrban} to the ultraweak space-time DPG discretization. The focus lies on what is specific to the present setting, namely the parameter dependence of the graph spaces and norms and the resulting uniform stability (\cref{sec:param}), the Kolmogorov $N$-width decay (\cref{sec:nwidth}), the role of the DPG optimal test spaces and the role of the built-in residual estimator in the reduced basis construction (\cref{sec:dpgrb,sec:remarks}).

\subsection{Parametrization and uniform well-posedness}
\label{sec:param}
We parametrize the potential and the right-hand side.
\begin{definition}
	\label{def:param}
	Let $\mathcal{P} \subseteq \R^P$ be a compact parameter set, $V_\mu \in L^\infty(\Omega_0;\R)$ be a parametrized potential, $A_\mu \coloneqq \ii \partial_t -\Delta_x + V_\mu$ the corresponding differential operator, and $b_\mu((u_\mu,q_\mu),v_\mu) \coloneqq \iO{u_\mu}{A_{h,\mu}^*v_\mu} + \duhm{q_\mu}{v_\mu}$ the corresponding sesquilinear form, where all spaces and norms from \cref{sec:setting} are defined with respect to $A_\mu$ and indexed by $\mu$ (e.g., $W_{h,\mu}^*$, $Q_\mu$). We are interested in solving the following problem.
	\begin{equation}
		\tag{\ensuremath{\mathrm{P}_\mu}} \label{morVarForm}
		\begin{aligned}
			\text{Given some } f_\mu \in (W_{h,\mu}^*)' \text{, find } (u_\mu, q_\mu) \in H \times Q_\mu \ \text{ s.t.}\\
			b_\mu((u_\mu,q_\mu),v_\mu) = f_\mu(v_\mu) \ \ \forall \, v_\mu \in W_{h,\mu}^*. \qquad \qquad
		\end{aligned}
	\end{equation}
\end{definition}
Concerning the parameter dependence, we impose the following standard assumptions.
\begin{assumption} \leavevmode
	\label{ass:affine}
	\begin{enumerate}[beginpenalty=10000]
		\item The $V_\mu$ and $f_\mu$ are parameter-separable, i.e., there exist functions $r_i:\mathcal{P} \to \R$ and $s_j:\mathcal{P} \to \C$, potentials $V_i \in L^\infty(\Omega_0;\R)$ and functionals $f_j \in (W_h^*)'$ with 
		\begin{align*}
			V_\mu = \sum_{i=1}^{K_V} r_i(\mu) V_i \ \ \text{ and } \ \ f_\mu = \sum_{j=1}^{K_f} s_j(\mu) f_j,
		\end{align*}
		where $K_V, K_f \in \N$ denote the number of affine components.
		\item The parameter dependencies $r_i: \mathcal{P} \to \R$ and $s_j: \mathcal{P} \to \C$ are continuous for all $1 \leq i \leq K_V$, $1 \leq j \leq K_f$.
	\end{enumerate}
\end{assumption}
In particular, the potentials are uniformly bounded with respect to the parameter, $C_V \coloneqq \sup_{\mu \in \mathcal{P}}\|V_\mu\|_\infty~<~\infty$. The affine parameter dependence forms the very foundation for an offline-online splitting, whereas continuity of the $r_i$ and $s_j$ enters the $N$-width result of \cref{sec:nwidth}.
\begin{notation}
	\label{notation:refnorms}
	In the following, let $X \coloneqq H \times Q$ and $Y \coloneqq W_h^*$ be equipped with the potential-\textit{free} reference norms, i.e., 
	\begin{align*}
		\nX{(u,q)}^2 
		&\coloneqq \nO{u}^2 + \inf_{w \in {\restr{D_h}{U}}^{-1} (\{q\})} \big(\nO{w}^2 + \nO{\ii \partial_t w - \Delta_x w}^2\big), \\
		\nY{v}^2 &\coloneqq \nO{v}^2 + \sum_{K \in \Omega_h} \nK{\ii \partial_t v - \Delta_x v}^2,
	\end{align*}
	so that all parameter-independent spaces and operators now denote the potential-free setting (e.g., $A = \ii\partial_t - \Delta_x$).
	In practice, these reference norms can be replaced by norms involving a suitably chosen mean potential in order to tighten the equivalence constants derived below, see \cref{sec:rom}.
\end{notation}
\begin{proposition}
	\label{propNormEquiv}
	Under \cref{ass:affine}, the trial spaces $X_\mu \coloneqq H \times Q_\mu$ and the test spaces $Y_\mu \coloneqq W_{h,\mu}^*$ coincide as sets with $X$ and $Y$, respectively, and their norms are uniformly equivalent,
	\begin{align}
		\label{eq:unifequiv}
		c_{X_\mu} \nX{x} \leq \nXm{x} \leq C_{X_\mu} \nX{x}, \qquad c_{X_\mu} \nY{v} \leq \nYm{v} \leq C_{X_\mu} \nY{v}
	\end{align}
	with $c_{X_\mu} \coloneqq (2 + \|V_\mu\|^2_\infty)^{-1/2} \geq (2 + C_V^2)^{-1/2}$ and $C_{X_\mu} \coloneqq (2 + \|V_\mu\|_\infty^2)^{1/2} \leq (2 + C_V^2)^{1/2}$. The same applies for the dual spaces $Y_\mu' = (W_{h,\mu}^*)'$.
\end{proposition}
\begin{proof}
	The graph spaces and boundary operators being independent of the potential (\cref{prop:potentialIndepGraphSpaces}) clearly implies the same for their parameter-dependent instances $Q_\mu$ and $W_{h,\mu}^*$. The equivalence for $\nYm{\cdot}$ is \eqref{eq:normequiv} with $V$ replaced by $V_\mu$, applied element by element and summed over $K \in \Omega_h$. The same holds for $\nWOm{\cdot}$ and consequently for $\nQm{\cdot}$ by monotonicity of the infimum. The dual norms inherit the equivalence by the definition of the operator norm. 
\end{proof}
Hence, \eqref{morVarForm} can be posed on the parameter-independent reference spaces $X$ and $Y$, which we do from now on. In the reference norms, the stability constants of \cref{thmwellposedweakSchr} depend on $C_V$, but \textit{uniformly} in $\mu$ and with explicit constants.
\begin{theorem}[Uniform well-posedness]
	\label{thm:unifwellposed}
	Suppose that \cref{ass:densityU,ass:affine} hold, and let $\beta_T$ be given by \eqref{eq:betaT}. Then, for every $\mu \in \mathcal{P}$ and all $(u,q) \in X$, $v \in Y$,
	\begin{align}
		\label{unifGammaBound}
		|b_\mu((u,q),v)| \leq (\sqrt{2} + \|V_\mu\|_\infty) \, \nX{(u,q)} \nY{v}
		&\leq (\sqrt{2} + C_V) \, \nX{(u,q)} \nY{v}, \\
		\label{unifInfSupConst}
		\adjustlimits\inf_{0 \neq x \in X} \sup_{0 \neq v \in Y} \frac{|b_\mu(x,v)|}{\nX{x} \nY{v}} \geq \beta_T \, (2+\|V_\mu\|^2_\infty)^{-1} &\geq \beta_T \, (2+C_V^2)^{-1}.
	\end{align}
	Consequently, \eqref{morVarForm} is well-posed for every $\mu \in \mathcal{P}$ with continuity and inf-sup constants depending only on $T$ and $C_V$, i.e., \textit{uniformly} in $\mu$.
\end{theorem}
\begin{proof}
	Fix $\mu \in \mathcal{P}$ and $x = (u,q)$. Since $b_\mu(x,v) = b(x,v) + \iO{u}{V_\mu v}$ with the potential-free form $b$ and $|\iO{u}{V_\mu v}| \leq \|V_\mu\|_\infty \nX{x} \nY{v}$, the continuity constant $\sqrt{2}$ of \cref{thmwellposedweakSchr} for $V = 0$ gives \eqref{unifGammaBound}. As \cref{ass:densityU} is independent of the potential, \cref{thmwellposedweakSchr} applied to $V_\mu$ yields, by \cref{rmk:constants}, an inf-sup constant $\beta_\mu \geq \beta_T$ with respect to $\nXm{\cdot}$ and $\nYm{\cdot}$. By \eqref{eq:unifequiv},
	\begin{align*}
		\sup_{0 \neq v \in Y} \frac{|b_\mu(x,v)|}{\nY{v}} \geq c_{X_\mu} \sup_{0 \neq v \in Y} \frac{|b_\mu(x,v)|}{\nYm{v}} \geq c_{X_\mu} \, \beta_T \, \nXm{x} \geq c_{X_\mu}^2 \, \beta_T \, \nX{x},
	\end{align*}
	and $c_{X_\mu}^{-2} = 2 + \|V_\mu\|_\infty^2$. Condition \eqref{N3} transfers from \cref{thmwellposedweakSchr} since the spaces coincide as sets.
\end{proof}
\begin{remark}[The role of the constants]
	\label{rmk:unifconst}
	Uniform inf-sup stability, a standard assumption in the reduced basis literature, is here a consequence of the potential-independent a~priori bound \eqref{eq:apriori}. All constants of the reduced basis scheme below depend on the parametrization only through
	\begin{align}
		\label{eq:refconst}
		\gamma_{\mathrm{ref}} \coloneqq \sqrt{2} + C_V \qquad \text{and} \qquad \beta_{\mathrm{ref}} \coloneqq \beta_T \, (2 + C_V^2)^{-1},
	\end{align}
	i.e., polynomially through the maximal magnitude $C_V$ of the potential.
\end{remark}

\subsection{Kolmogorov $N$-widths}
\label{sec:nwidth}
The approximability of the \textit{solution manifold} $\mathcal{M} \coloneqq \{x_\mu \ | \ \mu \in \mathcal{P}, \ x_\mu = (u_\mu, q_\mu) \text{ solves \eqref{morVarForm}}\} \subseteq X$ by $N$-dimensional subspaces is measured by the \textit{Kolmogorov $N$-widths}
\begin{align*}
	d_N(\mathcal{M}) \coloneqq \inf_{\substack{Z^N \subseteq X,\\\operatorname{dim}(Z^N) = N}} \ \sup_{\mu \in \mathcal{P}} \ \inf_{z^N \in Z^N} \nX{x_\mu - z^N}.
\end{align*}
The ingredients for the following \textit{exponential} decay have all been established above.
\begin{theorem}
	\label{thmNWidthConv}
	Under \cref{ass:densityU,ass:affine}, there exist constants $c,C > 0$ such that $d_N(\mathcal{M}) \leq C e^{-cN^{1/K_V}}$.
\end{theorem}
\begin{proof}
	This follows as in \cite[Theorem 1.3]{RBMUrban}, which requires parameter-independent reference spaces, uniform inf-sup stability in $\mu$, and an affine parameter dependence via continuous coefficient functions, all of which are assumed (\cref{ass:affine}) or have been verified (\cref{propNormEquiv,thm:unifwellposed}).
\end{proof}

This decay hinges on the way the parameter enters. Transport-dominated, wave-type and Schrödinger-like equations are known for slow $N$-width decay, as slow as $N^{-1/2}$ in some scenarios (see, e.g., \cite[Section 1.5.7]{RBMUrban}). It is also pointed out there that the underlying theorem does not apply to the reduced models of the dual formulations of the wave and Schrödinger equations \cite{ultrWeakWave, ultrWeakSchrodinger}, whose trial spaces are parameter-dependent already as sets, and that their rate of convergence is open. In contrast, the trial space $X = L^2(\Omega) \times D_h(U)$ of the present formulation is parameter-independent as a set and carries a uniformly equivalent topology, so that the theorem applies to the present parametrization of the potential and the right-hand side, at the price of the equivalence constants of \cref{propNormEquiv}.

\subsection{The DPG reduced basis method}
\label{sec:dpgrb}
The goal of the reduced basis method is to construct a low-dimensional, parameter-independent \textit{reduced space} $X^N \subseteq X$ that allows for very efficient solves of a \textit{reduced problem}. For a Petrov--Galerkin method, inf-sup stability of this problem requires suitable reduced test spaces. The principles are standard (cf.\ \cite{doGre, RBMUrban}), and it is the DPG framework that renders them computationally accessible, see \cref{sec:remarks}.
\paragraph{Truth discretization}
As the ideal DPG method inherits stability from the continuous level, it remains to consider the \textit{practical} method for \eqref{morVarForm}. We fix a trial space $X_\scn \subseteq X$ and an enriched test space $Y^{\scDelta}_\scn \subseteq Y$, both independent of $\mu$, and define the parameter-dependent trial-to-test map analogous to \eqref{Tdiscrdef} by
\begin{align}
	\label{mapTscnmu}
	T_{\scn,\mu}: X_\scn \to Y^\scDelta_\scn, \quad 
	\iY{T_{\scn,\mu} x_\scn}{y^\scDelta_\scn} = b_\mu(x_\scn,y^\scDelta_\scn) 
	\quad \text{for all } y^\scDelta_\scn \in Y^{\scDelta}_\scn.
\end{align}

\begin{proposition}
	\label{propPracticalMu}
	Assume that, for every $\mu \in \mathcal{P}$, there exists a Fortin operator $\Pi_\mu: Y \to Y^{\scDelta}_\scn$ with $b_\mu(x_\scn, y - \Pi_\mu y) = 0$ for all $x_\scn \in X_\scn$ and $y \in Y$, and that $\nu \coloneqq \sup_{\mu \in \mathcal{P}} \|\Pi_\mu\| < \infty$. Then, the \textit{truth solution} $x_{\scn,\mu} \in X_\scn$ of the practical DPG method for \eqref{morVarForm} with test space $T_{\scn,\mu}(X_\scn)$ satisfies the uniform quasi-optimality condition
	\begin{equation}
		\label{quasiOptMu}
		\nX{x_\mu - x_{\scn,\mu}} 
		\leq \frac{\nu \, \gamma_{\mathrm{ref}}}{\beta_{\mathrm{ref}}}
		\, \inf_{z_\scn \in X_\scn} \nX{x_\mu - z_\scn}.
	\end{equation}
\end{proposition}
\begin{proof}
	This is \cref{thmwellposedPracticalSchr}, applied for every $\mu \in \mathcal{P}$, with the constants $\gamma_\mathrm{ref}$ and $\beta_\mathrm{ref}$ from \eqref{eq:refconst}.
\end{proof}
Hence, the parametrization causes no additional difficulty for the practical method beyond a uniform bound on the Fortin operators for the individual potentials, whose existence is open already for the free equation (cf.\ \cref{rmkDiscreteStab}).
\paragraph{Reduced problem and error estimation}
Given a reduced trial space $X^N \subseteq X_\scn$, we define the parameter-dependent reduced test space $Y_\mu^N \coloneqq T_{\scn,\mu}(X^N)$ via \eqref{mapTscnmu} and consider the following \textit{reduced problem}.
\begin{equation}
	\tag{\ensuremath{\mathrm{R}_\mu}} \label{rbVarForm}
	\begin{aligned}
		\text{Given some } f_\mu \in Y' \text{, find } (u_\mu^N, q_\mu^N) \in X^N \ \text{ s.t.}\\
		b_\mu((u^N_\mu,q^N_\mu),v^N) = f_\mu(v^N) \ \ \forall\, v^N \in Y^N_\mu. \qquad 
	\end{aligned}
\end{equation}

\begin{proposition}[Reduced normal equations]
	\label{propNormalEq}
	Let $X^N \subseteq X_\scn$ be spanned by $x^N_1, \dots, x^N_N$. Then, \eqref{rbVarForm} can be expressed in matrix form as
	\begin{align}
		\label{rbMatrixForm}
		(\boldsymbol{\mathcal{B}}_\mu^N)^H \boldsymbol{\mathcal{G}}^{-1} \boldsymbol{\mathcal{B}}_\mu^N \,
		\boldsymbol{x}_\mu^N
		= 
		(\boldsymbol{\mathcal{B}}^N_\mu)^H \boldsymbol{\mathcal{G}}^{-1} \boldsymbol{\mathcal{F}}_\mu,
	\end{align}	
	where $\boldsymbol{\mathcal{G}} \in \C^{\mathcal{N}^{\scDelta} \times \mathcal{N}^{\scDelta}}$, $\mathcal{N}^\scDelta \coloneqq \operatorname{dim}(Y^\scDelta_\scn)$, denotes the Gramian of a basis $(y_{\scn,i})_i$ of $Y_\scn^{\scDelta}$ with respect to the reference inner product, the matrix $\boldsymbol{\mathcal{B}}_\mu^N \in \C^{\mathcal{N}^{\scDelta} \times N}$ has the entries $b_\mu(x^N_j, y_{\scn,i})$, the vector $\boldsymbol{\mathcal{F}}_\mu$ has the entries $f_\mu (y_{\scn,i})$, and $\boldsymbol{x}_\mu^N$ is the coefficient vector of $x_\mu^N$ with respect to $(x^N_i)_i$. Under the assumptions of \cref{propPracticalMu}, $x_\mu^N$ satisfies
	\begin{equation}
		\label{quasiOptReduced}
		\nX{x_\mu - x_\mu^N} 
		\leq \frac{\nu \, \gamma_{\mathrm{ref}}}{\beta_{\mathrm{ref}}} \, \inf_{z^N \in X^N} \nX{x_\mu - z^N}.
	\end{equation}
\end{proposition}
\begin{proof}
	The matrix representation follows since $Y_\mu^N$ is exactly represented by the columns of $\boldsymbol{\mathcal{G}}^{-1} \boldsymbol{\mathcal{B}}_\mu^N$. As $Y_\mu^N$ is the span of the supremizers in $Y_\scn^{\scDelta}$, i.e., $\nY{T_{\scn,\mu} z^N} = \sup_{y \in Y^N_\mu} |b_\mu(z^N, y)| / \nY{y}$ for $z^N \in X^N$, the discrete inf-sup constant $\beta_\mathrm{ref}/\nu$ of \cref{thmwellposedPracticalSchr} is inherited on $X^N \times Y^N_\mu$, and \eqref{quasiOptReduced} follows as in \cref{propPracticalMu}.
\end{proof}
\begin{definition}
	\label{def:estimator}
	For any $\mu \in \mathcal{P}$ and subspace $X^N \subseteq X_\scn$ with reduced solution $x_\mu^N \in X^N$ of \eqref{rbVarForm}, we define the estimator
	\begin{align}
		\label{greedyErrorEstimator}
		\eta_\scn(\mu, X^N) \coloneqq \nYDp{B_\mu x_\mu^N - f_\mu} = \sup_{0 \neq y^\scDelta_\scn \in Y_\scn^\scDelta} \frac{|b_\mu(x_\mu^N,y^\scDelta_\scn) - f_\mu(y^\scDelta_\scn)|}{\nY{y^\scDelta_\scn}}. \!\!\!
	\end{align}
\end{definition}
\begin{proposition}
	\label{propQuasiTight}
	Under the assumptions of \cref{propPracticalMu}, the estimator $\eta_\scn$ bounds the error to the \textit{exact} solution, namely
	\begin{align*}
		\frac{1}{\gamma_{\mathrm{ref}}} \, \eta_\scn(\mu,X^N) \leq \nX{x_\mu - x_\mu^N} \leq \frac{1}{\beta_{\mathrm{ref}}} \, \bigl(\nu \, \eta_\scn(\mu,X^N) + \operatorname{osc}(f_\mu)\bigr)
	\end{align*}
	with the data oscillation $\operatorname{osc}(f_\mu)$ from \eqref{oscf} for $\Pi = \Pi_\mu$, and combined with \eqref{quasiOptReduced} the \textit{quasi-tightness} condition
	\begin{equation*}
		\frac{\beta_{\mathrm{ref}}}{\nu \, \gamma_{\mathrm{ref}}^2} \, \eta_\scn(\mu,X^N) \leq \inf_{z^N \in X^N} \nX{x_\mu - z^N} \leq \frac{1}{\beta_{\mathrm{ref}}} \, \bigl(\nu \, \eta_\scn(\mu,X^N) + \operatorname{osc}(f_\mu)\bigr).
	\end{equation*}
\end{proposition}
\begin{proof}
	The first statement is \eqref{apostApproxReliab} for $\tilde{x}_\scn = x_\mu^N$ with the uniform constants from \cref{thm:unifwellposed}. The lower bound of the second follows from \eqref{quasiOptReduced}, the upper bound from $\inf_{z^N \in X^N} \nX{x_\mu - z^N} \leq \nX{x_\mu - x_\mu^N}$. 
	\end{proof}
Ignoring the data oscillation term, \cref{propQuasiTight} shows that $\eta_\scn$ is \textit{tight} in the sense of \cite{doGre}, i.e., equivalent to the best approximation error $\inf_{z^N \in X^N} \nX{x_\mu - z^N}$ uniformly in $\mu$ and $X^N$.
\paragraph{Offline-online decomposition}
The reduced problem is solved online in the form \eqref{rbMatrixForm}, as the mixed system \eqref{eq:blocksystem} over $X^N$ retains the full-order error representation function. Both \eqref{rbMatrixForm} and the estimator \eqref{greedyErrorEstimator} admit the standard offline-online decomposition \cite{RBMUrban}, since, by the affine dependence of $b_\mu$ and $f_\mu$ and the parameter-independent Gramian, the columns of $\boldsymbol{\mathcal{G}}^{-1} \boldsymbol{\mathcal{B}}_\mu^N$ and the Riesz representative of the residual, whose norm is $\eta_\scn(\mu,X^N)$ (cf.\ \cref{propMixedResidPract}), are $\mu$-weighted combinations of precomputable elements of $Y^\scDelta_\scn$. This is the reason for defining $T_{\scn,\mu}$ in \eqref{mapTscnmu} via the reference inner product, the inverse of a $\mu$-dependent Gramian $\boldsymbol{\mathcal{G}}_\mu$ admitting no affine decomposition.
\paragraph{Greedy algorithm}
The reduced spaces are generated by the standard greedy algorithm \cite{RBMUrban, wGreConv}. Given a finite training set $\mathcal{S} \subseteq \mathcal{P}$ and starting from $X^0 \coloneqq \{0\}$, the $N$-th step selects $\mu_N \coloneqq \arg\max_{\mu \in \mathcal{S}} \eta_\scn(\mu, X^{N-1})$, computes the snapshot $x_{\scn,\mu_N}$ by the practical DPG method of \cref{propPracticalMu}, and sets $X^N \coloneqq X^{N-1} + \operatorname{span}\{x_{\scn,\mu_N}\}$, until $\max_{\mu \in \mathcal{S}} \eta_\scn(\mu, X^N)$ drops below a prescribed tolerance. By the tightness of $\eta_\scn$, the greedy spaces inherit the exponential decay of \cref{thmNWidthConv} with the degraded exponent $1/(1+K_V)$ by \cite[Theorem 3.5]{wGreConv}, a guarantee that holds only down to the level of the truth error and the data oscillation $\operatorname{osc}(f_\mu)$.

\subsection{Discussion}
\label{sec:remarks}
As in \cite{minResMixedRB, LSFEMRB}, the estimator $\eta_\scn$ measures the error to the \textit{infinite-dimensional} solution. It is rooted in the dual norm of the \textit{enriched} test space $Y_\scn^\scDelta$, so that $\nYDp{B_\mu x_{\scn,\mu} - f_\mu} \neq 0$ even for the truth solutions. The estimator can only be driven down to this truth baseline, whereas the errors $\nX{x_\mu^N - x_{\scn,\mu}}$ between reduced and truth solution may drop significantly below. As a DPG-specific alternative to a fixed tolerance, the truth residual $\nY{e_{\scn, \mu}^\scDelta}$ obtained with every snapshot is a lower bound for the estimator, so that termination could be coupled to their difference.

Supremizer-based schemes are known to be computationally demanding for standard continuous formulations, requiring global solves for all trial basis vectors, whereas the discontinuous test space of DPG methods renders the Gramian block-diagonal and the supremizer computation element-wise. Together with static condensation, which reduces the global full-order solve to the interface unknowns (\cref{rmkStaticCond}), this makes supremizer-based reduced basis schemes accessible for large-scale problems.

Closely connected is the dual approach of \cite{brunkenTranspEq, MORUltraweakTransport, RBMUrban}, which obtains a trial space perfectly matching the test space as $X = B^*(Y)$, requires no Riesz map since $X = L^2(\Omega)$, and constructs the reduced \textit{test} space greedily without any supremizer computations. Here, however, the interface space $Q = D_h(U)$ carries the quotient norm $\nQ{q} = \nWhsp{q}$ of the broken second-order graph space, so that evaluating $B^*$ requires the \textit{same} Gramian inversion as the test space construction and the dual scheme yields no advantage. Moreover, the parameter dependence of the spaces in \cite{RBMUrban} renders their residual not offline-online decomposable \cite[Section 1.5.4]{RBMUrban}, which requires additional care.
\section{Numerical experiments}
\label{sec:numerics}
All numerical experiments employ NGSolve \cite{NGSolve} for the discretization and the full-order solves, embedded in pyMOR \cite{pyMORPaper} to drive the greedy construction. The full-order systems are solved via static condensation of \eqref{eq:blocksystem} (\cref{rmkStaticCond}) with the direct sparse solver \textit{PARDISO} for the Schur complement, which provides the built-in error estimator $\nY{e^\scDelta_\scn}$ at no additional cost. The interface unknowns are discretized as in \cref{sec:discrete1D}, the fluxes $q^{\vtick}$ by facet polynomials of degree $p$ (\textit{FacetFESpace}), whereas \cite{DPGSchr} employed normal traces of $H(\mathrm{div})$-conforming elements of degree $p-1$, one degree lower than required by \cref{thmerror1DSchrDPG}. We evaluate both flux degrees. The complete source code is publicly available at \cite{mycode}.
\subsection{Full-order model}
\label{sec:fom}
\paragraph{Validation without potential}
We first reproduce Case (a) of \cite{DPGSchr}, a dispersive pulse in an optical fiber modeled by the complex Gaussian
\begin{equation}
	u_{(a)}(x,t) = \frac{M T_0}{\sqrt{T_0^2 - \ii \beta t}} \exp\left( - \frac{x^2}{T_0^2 - \ii \beta t} \right), \qquad M = T_0 = 1.5, \ \beta = 2.5,
	\label{case_a_exact}
\end{equation}
on $\Omega = (0,1)^2$, where $\partial_{xx} u$ is scaled by $\frac{\beta}{2}$, which does not affect the preceding analysis, and the right-hand side $f = \ii \partial_t u - \frac{\beta}{2} \partial_{xx} u$ is derived analytically. The inhomogeneous boundary and initial data on $\Gamma$ are carried by $Q_\scn^{\tickplus}$ and imposed by the usual modification of the right-hand side of \eqref{eq:blocksystem}. The second example of \cite{DPGSchr}, Case (b), is stated as a real-valued expression \cite[(5.2)]{DPGSchr}, whereas the corresponding plots display radial phase oscillations, so that an exact comparison is not possible.

\begin{figure}[htbp]
	\centering
	\includegraphics[width=0.48\textwidth]{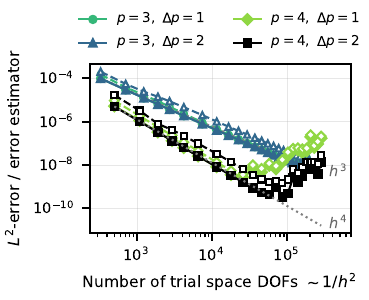}
	\hfill
	\includegraphics[width=0.48\textwidth]{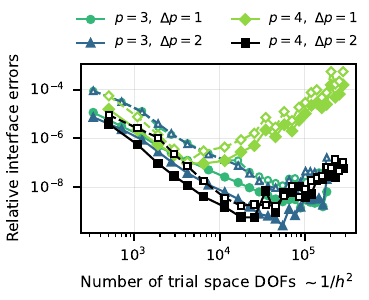}
	\caption{Case (a) with the full flux degree $p$. Volume $L^2$-errors (solid) and estimator (dashed) on the left, relative interface errors of $q^{\tickplus}$ (solid) and $q^{\vtick}$ (dashed) on the right}
	\label{fig:case_a}
\end{figure}

\Cref{fig:case_a} displays the $L^2$-errors of $u_\scn$, the relative skeleton $L^2$-errors of $q^{\tickplus}$ and $q^{\vtick}$, and the built-in estimator $\nY{e_\scn^\scDelta}$ for $p \in \{3,4\}$, $\Delta p \in \{1,2\}$, and the full flux degree $p$ required by \cref{thmerror1DSchrDPG}. By \eqref{apostApproxReliab}, the estimator is no strict bound for the volume error. The skeleton $L^2$-norms merely serve as computable proxies for the intractable quotient norm $\nQ{\cdot}$ to verify that the interface variables converge. With the reduced flux degree $p-1$ of \cite{DPGSchr}, whose plots we do not display, we reproduce the results reported there, and two observations are in order. First, as already noted in \cite{DPGSchr}, the observed rates exceed the $\mathcal{O}(h^{p-1})$ predicted by \cref{thmerror1DSchrDPG} by a full polynomial degree, where for $p = 4$ we even observe a clean $\mathcal{O}(h^4)$ rate instead of a rate between $\mathcal{O}(h^3)$ and $\mathcal{O}(h^4)$. Second, the error curves eventually stagnate and fluctuate due to ill-conditioning \cite[Section~5]{DPGSchr}, but only at errors about two orders of magnitude below those reported there, which we attribute to static condensation and the direct solver. Both observations carry over to the full flux degree in \cref{fig:case_a}. The main difference concerns stability. The configuration $p = 4$, $\Delta p = 1$, which does not display any stability issues for the reduced flux degree, becomes unstable for the full one, clearly visible in the interface errors, which fluctuate without any consistent asymptotic rate. The larger trial space demands more from the enriched test space $W_\scn^{*\scDelta}$, and $\Delta p = 1$ is simply insufficient, a stability issue not visible in \cite{DPGSchr}. We therefore use $\Delta p = 2$ in \cref{sec:rom}.

\paragraph{Harmonic oscillator potential}
We now add a potential and consider as Case (c) the harmonic oscillator $V(x) = \frac{1}{2} m \omega^2 (x - x_c)^2$, an important model in solid-state physics and quantum field theory \cite[Section 11.1]{quantMath}. The exact solution of the homogeneous equation, $f = 0$, is the Gaussian wave packet of \cite[Problem 2.49]{introductQuantMech}. In the units of \cref{Schrödingerclass}, i.e., $\hbar = 1$ and $m = \frac{1}{2}$ in \cite[(1.1)]{introductQuantMech}, and with $\xi \coloneqq x - x_c$, the potential reads $V(x) = \frac{\omega^2}{4} \xi^2$ and
\begin{equation}
	u_{(c)}(x,t) = \left(\frac{\omega}{2\pi}\right)^{\!1/4} \exp\left[ -\frac{\omega}{4} \left( \xi^2 + \frac{x_0^2}{2} \left(1 + e^{2\ii\omega t}\right) - 2\ii t - 2 x_0 \xi e^{\ii\omega t} \right) \right],
	\label{case_c_exact}
\end{equation}
solves $\ii \partial_t u - \partial_{xx} u + Vu = 0$ for any $x_0 \in \R$, where we have taken the complex conjugate of the textbook expression, as \cref{Schrödingerclass} is the complex conjugate of \cite[(1.1)]{introductQuantMech}. We choose $x_c = \frac{1}{2}$, $x_0 = 0.2$ and $\omega = 8\pi$ on $\Omega = (0,1)^2$, i.e., four oscillations within the time horizon, and impose the data of $u_{(c)}$ on $\Gamma$ as in Case (a).

\begin{figure}[htbp]
	\centering
	\begin{subfigure}[b]{0.48\textwidth}
		\centering
		\includegraphics[width=\textwidth]{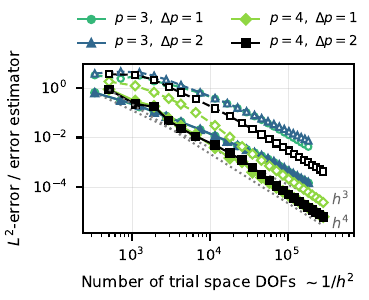}
		\caption{Volume $L^2$-errors (solid) and error estimator (dashed)}
	\end{subfigure}
	\hfill
	\begin{subfigure}[b]{0.48\textwidth}
		\centering
		\includegraphics[width=\textwidth]{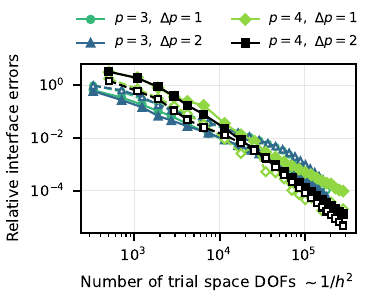}
		\caption{Relative interface variable errors of $q^{\tickplus}$ (solid) and $q^{\vtick}$ (dashed)}
	\end{subfigure}
	\caption{Convergence results for Case (c), the harmonic oscillator potential with $\omega = 8\pi$}
	\label{fig:case_c}
\end{figure}

\Cref{fig:case_c} displays the same quantities as \cref{fig:case_a}. We again observe rates of about $\mathcal{O}(h^p)$ on the finer meshes with practically identical volume errors for $\Delta p = 1$ and $\Delta p = 2$. Hence, the potential does not affect the convergence order, which again exceeds the guarantee of \cref{thmerror1DSchrDPG}.

\subsection{Reduced-order model}
\label{sec:rom}
We finally demonstrate the explicit combination of the DPG method with the greedy reduced basis scheme, implemented as a ``hybrid'' approach, i.e., we utilize static condensation of \eqref{eq:blocksystem} for computing the high-dimensional truth solutions, while driving the greedy algorithm via the normal equations \eqref{rbMatrixForm}, which makes the supremizer-based reduction accessible at this scale, see \cref{sec:remarks}. For the full-order model, we employ a uniform $90 \times 90$ mesh of square elements ($h = 1/90$) on $\Omega = (0,1)^2$, the trial space polynomial degree $p=3$ with the full flux degree for $Q_\scn^{\vtick}$, and $\Delta p = 2$, i.e., $146\,701$ trial degrees of freedom. The greedy algorithm is executed over a randomly sampled training set of $|\mathcal{S}| = 10^4$ parameters until the estimator values $\{\eta_\scn(\mu,X^N)\}_\mu$ uniformly drop below $\epsilon_{\text{tol}} = 10^{-4}$ or $100$ iterations are reached. As mentioned in \cref{notation:refnorms}, the Gramian $\boldsymbol{\mathcal{G}}$ and the Gramian block of \eqref{eq:blocksystem} are evaluated with respect to the mean-parameter potential, so that the truth solutions are the practical DPG solutions of \cref{propPracticalMu} and the truth residual $\nY{e^\scDelta_{\scn,\mu}}$ is a lower bound for $\eta_\scn(\mu, X^N)$. All errors below are measured in the $L^2$-norms of the respective components, the interface variables on the skeleton, and we consider two cases with homogeneous data on $\Gamma$.
\begin{description}
	\item[Case 1 ($K_V = 2$).] The parameter vector is $\mu \in \mathcal{P} = [0, 150]^2$,
	\begin{align*}
		V_\mu(x) &= \mu_1 \exp\left(-200(x - 0.3)^2\right) + \mu_2 \exp\left(-200(x - 0.7)^2\right), \\
		f(x,t) &= \sin(\pi x) \sin(10 t).
	\end{align*}
	\item[Case 2 ($K_V = 3$).] The parameter vector is $\mu \in \mathcal{P} = [0, 200] \times [0, 150] \times [-40, 40]$,
	\begin{align*}
		V_\mu(x) &= \mu_1(x - 0.5)^2 + \mu_2 \exp\left(-100(x - 0.5)^2\right) + \mu_3(x - 0.5), \\
		f(x,t) &= \exp\left(-150(x - 0.2 - 0.6t)^2\right) \sin(\pi t).
	\end{align*}
\end{description}
\begin{figure}[htbp]
	\centering
	\begin{subfigure}[b]{0.48\textwidth}
		\centering
		\includegraphics[width=\textwidth]{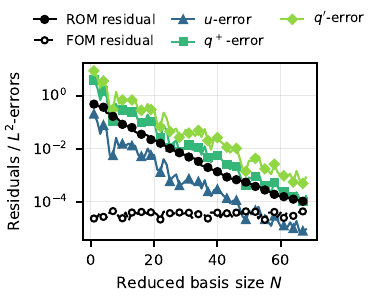}
		\caption{Case 1, $N = 68$}
	\end{subfigure}
	\hfill
	\begin{subfigure}[b]{0.48\textwidth}
		\centering
		\includegraphics[width=\textwidth]{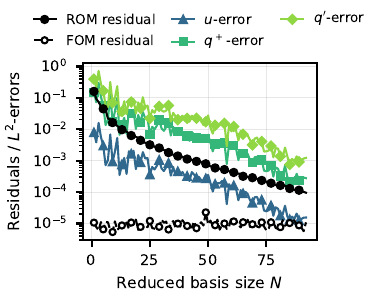}
		\caption{Case 2, $N = 92$}
	\end{subfigure}
	\caption{Greedy basis construction for Cases 1 and 2}
	\label{fig:spacetime_greedy_convergence}
\end{figure}
\Cref{fig:spacetime_greedy_convergence} displays, along the greedy construction, the maximal estimator over the training set, the truth residual of the selected snapshot and the $L^2$-errors between reduced and truth solution at the selected parameter. The estimator decays exponentially once an initial exploration of the parameter domain is concluded, at a slower rate in Case 2, and in both cases in agreement with the rate $e^{-cN^{1/K_V}}$ of \cref{thmNWidthConv}. The volume and the interface errors decay at comparable rates, and the $L^2$-error of $u$ stays a factor between $2$ and $60$ below the estimator, although the worst-case inf-sup constant in \cref{thm:unifwellposed} with the mean-parameter reference norm is only of order $(C_V/2)^{-2} \approx 2 \cdot 10^{-4}$. The estimator approaches the truth residual, the baseline of \cref{sec:remarks}, from above, whereas the $L^2$-errors of $u$ drop below the baseline in Case 1, so that $\epsilon_{\text{tol}}$ has to exceed the truth residuals. The condition numbers of the reduced normal equations increased at most to $3 \cdot 10^4$, thus not causing any conditioning issues.

\begin{table}[htbp]
	\centering
	\caption{Performance of the reduced-order models over $N_{\text{test}} = 200$ random test parameters per case (wall-clock solve times, relative errors in the respective $L^2$-norms)}
	\label{tab:spacetime_rom_performance}
	\footnotesize
		\begin{tabular}{l c c c c c c}
			\toprule
			& \multicolumn{3}{c}{\textbf{Case 1} ($N=68$)} & \multicolumn{3}{c}{\textbf{Case 2} ($N=92$)} \\
			\cmidrule(lr){2-4} \cmidrule(lr){5-7}
			\textbf{Performance indicator} & \textbf{Min.} & \textbf{Median} & \textbf{Max.} & \textbf{Min.} & \textbf{Median} & \textbf{Max.} \\
			\midrule
			ROM solve time [ms] & $0.69$ & $0.82$ & $1.61$ & $1.12$ & $1.48$ & $6.86$ \\
			Computational speedup & $9082\times$ & $17775\times$ & $21141\times$ & $2119\times$ & $9790\times$ & $12910\times$ \\
			Relative $u$-error [\%] & $0.0004$ & $0.0122$ & $0.0603$ & $0.0200$ & $0.1112$ & $0.3710$ \\
			Relative $q^{\tickplus}$-error [\%] & $0.0004$ & $0.0093$ & $0.0492$ & $0.0098$ & $0.0614$ & $0.3340$ \\
			Relative $q^{\vtick}$-error [\%] & $0.0007$ & $0.0174$ & $0.0861$ & $0.0189$ & $0.0931$ & $0.4003$ \\
			ROM resid. / FOM resid. & $1.01$ & $1.44$ & $2.69$ & $1.46$ & $4.66$ & $12.45$ \\
			\bottomrule
		\end{tabular}
\end{table}

\Cref{tab:spacetime_rom_performance} evaluates the final reduced-order models of dimension $N = 68$ and $N = 92$ on $N_\text{test} = 200$ random parameters, where the reduced solutions reproduce the truth solutions with relative errors below $0.1\,\%$ and $0.5\,\%$ at solve times of about one millisecond against full-order solve times between $14.4$ and $15.7$ seconds for the condensed system. Moreover, the ratio between ROM and FOM residual is at least one, as it must be, and larger for Case 2, which terminated further above the baseline.

\begin{figure}[htbp]
	\centering
	\includegraphics[width=0.48\textwidth]{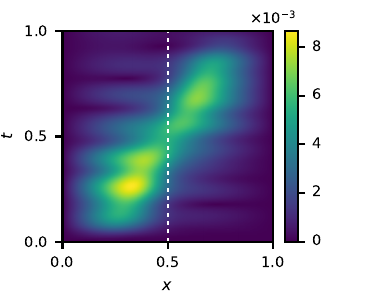}
	\hfill
	\includegraphics[width=0.48\textwidth]{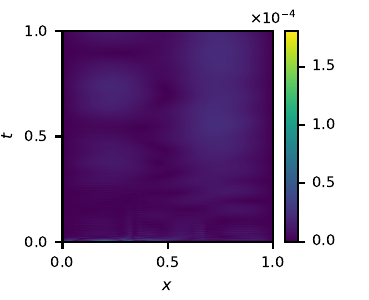}
	\caption{Case 2 at the test parameter with the largest relative $u$-error. Absolute value of the truth solution (left, $x = 0.5$ dashed) and its difference to the reduced one (right)}
	\label{fig:spacetime_solution}
\end{figure}

\Cref{fig:spacetime_solution} displays the truth solution of Case 2 at the test parameter with the largest relative $u$-error, which follows the moving source, together with its pointwise difference to the reduced solution.

\section{Conclusion}
\label{sec:conclusion}
We extended the ultraweak space-time DPG method of \cite{DPGSchr} to bounded, time-independent potentials, which enter the well-posedness proof only through the spectral decomposition of the spatial operator and leave the graph spaces and the density assumption unchanged. The resulting inf-sup constant in the graph norms is uniform in the potential, which makes the parametrized problem accessible to the reduced basis theory on parameter-independent reference spaces and yields exponentially decaying Kolmogorov $N$-widths. The optimal test spaces and the built-in residual of the DPG method provide a feasible supremizer-based reduced test space and an error estimator with respect to the infinite-dimensional solution. In the experiments, the flux degree of the a~priori estimate required the enrichment $\Delta p = 2$ for $p = 4$, the full-order method converged at the rates of the free equation also with a potential, one order above this estimate, and reduced models of dimension $68$ and $92$ reproduced the truth solutions with relative errors below $0.5\,\%$ at speedups between $2 \cdot 10^3$ and $2 \cdot 10^4$, with an estimator decaying exponentially in the reduced dimension.

Several questions remain open. Time-dependent potentials, permitted for instance by the framework of \cite{ultrWeakSchrodinger}, require well-posedness arguments beyond the presented spatial spectral decomposition. Further, as for the free equation, a Fortin operator for the broken second-order graph norm is missing, so that the observed robustness of the practical method with respect to the potential lacks a theoretical foundation. Finally, several spatial dimensions, more general parametrizations, which affect the graph spaces as sets, and the application of the DPG-RBM combination to other parametric problems are natural continuations of this work.

\appendix
\section{The discontinuous Petrov--Galerkin method}
\label{sec:dpg}
The discontinuous Petrov--Galerkin (DPG) method was introduced by Demkowicz and Gopalakrishnan in 2010 \cite{DPGpart1}. Its defining features are a \textit{discontinuous} test space, which decouples the computation of \textit{optimal} test functions into local, element-wise problems, and, since these optimal test functions are not computable exactly, a \textit{practical} variant based on a sufficiently enriched discrete test space. We collect the results needed in \cref{sec:schroedinger,sec:mor,sec:numerics}, following the survey \cite{DPGActaNumerica} and references therein.
\begin{notation}
	\label{notation:dpg}
	Let $(X,\iX{\cdot}{\cdot})$ and $(Y,\iY{\cdot}{\cdot})$ be complex Hilbert spaces, where $Y$ consists of functions $\Omega \to \C$ on the domain $\Omega \subseteq \R^d$ with mesh $\Omega_h$ as in \cref{sec:schroedinger}, and let $\dY{\cdot}{\cdot}$ denote the dual pairing on $Y' \times Y$.
	Subsequently, $Y= \prod_{K \in \Omega_h} Y(K)$ consists of local Hilbert spaces on the single mesh elements and is equipped with the product topology, $\nY{\cdot}^2 = \sum_{K \in \Omega_h}\nYK{\cdot}^2$, so that the functions in $Y$ do not fulfill any global continuity constraint. $X_\scn \subseteq X$ denotes some finite-dimensional subspace.
	Finally, $b: X \times Y \to \C$ denotes a continuous sesquilinear form, antilinear in its second variable, and $B: X \to Y'$ the corresponding bounded linear operator, $\dY{Bx}{y} \coloneqq b(x,y)$. 
\end{notation}

We are interested in solving the following variational formulation. 
\begin{equation}
	\begin{aligned}
		\text{Given some } f \in Y' \text{, find } x \in X \ \text{ s.t.}\\
		b(x,y) = f(y) \ \ \forall \, y \in Y. \ \ \ \ \ \
	\end{aligned}
	\tag{\text{P}} \label{varForm}
\end{equation}
We assume that \eqref{varForm} satisfies the prerequisites of the \textit{Ne\v{c}as theorem} (see, e.g., \cite[Theorem 1.1]{RBMUrban}), that is
\begin{align}
	\label{N1} \tag{N1} &\exists \, \gamma < \infty \textnormal{ s.t. } |b(x,y)| \leq \gamma \nX{x} \nY{y} \ \ \forall \, x \in X, y \in Y, \\
	\label{N2} \tag{N2} &\exists \, \beta > 0 \textnormal{ s.t. } \infsupXY{x}{y}{X}{Y} \geq \beta, \\
	\label{N3} \tag{N3} &\forall \, 0 \neq y \in Y \ \ \exists x \in X  \text{ s.t. } b(x,y) \neq 0.
\end{align}
Then, $B: X \to Y'$ is an isomorphism with $\|B\| \leq \gamma$ and $\|B^{-1}\| \leq \frac{1}{\beta}$, and we refer to $\gamma$ and $\beta$ as the continuity and inf-sup constants of $b$ on $X \times Y$.

\paragraph{The ideal DPG method}
The \textit{ideal} DPG method for \eqref{varForm} reads as follows.
\begin{equation}
	\begin{aligned}
		\text{Given some } f \in Y' \text{, find } x_\scn \in X_\scn \ \text{ s.t.}\\
		b(x_\scn,y_\scn) = f(y_\scn) \ \ \forall\, y_\scn \in T(X_\scn), \ \ \ \ \ \
	\end{aligned}
	\tag{\text{I-DPG}} \label{DPG}
\end{equation}
Here, $T: X \to Y$ maps $x$ to the \textit{Riesz representative} of $Bx \in Y'$, i.e.,
\begin{align}
	\label{Tdef}
	\iY{Tx}{y} = b(x,y) \text{ for all } y \in Y.
\end{align}
Since $\sup_{y \in Y} |b(x,y)|/\nY{y} = \nY{Tx}$ for all $x \in X$, the \textit{optimal} test space $T(X_\scn)$ consists of \textit{supremizers}, and inf-sup stability on $X \times Y$ is inherited by $X_\scn \times T(X_\scn)$ with the same constant. This yields the following result (see \cite[Theorem 2.4]{DPGActaNumerica}).
\begin{theorem}
	\label{thmideal}
	Let $x \in X$ denote the unique solution of \eqref{varForm}. Then, \eqref{DPG} is uniquely solved by some $x_\scn \in X_\scn$ satisfying the quasi-optimality condition 
	\begin{align}
		\label{bestapprox}
		\nX{x - x_\scn} \leq \frac{\gamma}{\beta} \inf_{z_\scn \in X_\scn} \nX{x - z_\scn}.
	\end{align}
\end{theorem}
\paragraph{The practical DPG method and a posteriori error control}
Computing $T(X_\scn)$ requires the Riesz representatives \eqref{Tdef} in the infinite-dimensional space $Y$ and is thus intractable. Hence, we choose a finite-dimensional \textit{enriched} test space $Y^{\scDelta}_\scn \subseteq Y$, in practice broken piecewise polynomials of degree $p + \Delta p$ with $p$ the degree of $X_\scn$ and an enrichment $\Delta p \geq 1$, and consider $T_\scn: X_\scn \to Y^{\scDelta}_\scn$ satisfying
\begin{align}
	\label{Tdiscrdef}
	\iY{T_\scn x_\scn}{y^\scDelta_\scn} = b(x_\scn,y^\scDelta_\scn) \ \text{ for all } y^\scDelta_\scn \in Y^{\scDelta}_\scn,
\end{align}	
i.e., the $Y$-orthogonal projection of the supremizers onto $Y^{\scDelta}_\scn$. With the \textit{approximate optimal test space} $Y_\scn \coloneqq T_\scn(X_\scn)$, the \textit{practical} DPG method reads as follows.
\begin{equation}
	\begin{aligned}
		\text{Given some } f \in Y' \text{, find } x_\scn \in X_\scn \ \text{ s.t.}\\
		b(x_\scn,y_\scn) = f(y_\scn) \ \ \forall \, y_\scn \in Y_\scn. \ \ \ \
	\end{aligned}
	\tag{\text{P-DPG}} \label{approxidealDPG}
\end{equation}

\begin{theorem}[{\cite[Theorem 2.1 and its proof]{practicalDPG}}]
	\label{thmwellposedPracticalSchr}
	Assume that there exists a linear \textit{Fortin} operator $\Pi: Y \to Y_\scn^\scDelta$ such that 
	\begin{equation}
		\label{Fortin}
		\nu \coloneqq \|\Pi\| < \infty \qquad \text{and} \qquad b(x_\scn,y - \Pi y) = 0 \quad \forall \, x_\scn \in X_\scn, \, y \in Y.
	\end{equation}
	Then, \eqref{approxidealDPG} fulfills the discrete inf-sup condition
	\begin{align*}
		\infsupXY{x_\scn}{y_\scn}{X_\scn}{Y_\scn} \geq \beta_\scn \coloneqq \frac{\beta}{\nu},
	\end{align*}and therefore admits a unique solution $x_\scn \in X_\scn$ satisfying
	\begin{align}
		\label{quasiOptPractical}
		\nX{x - x_\scn} \leq \frac{\nu \gamma}{\beta} \inf_{z_\scn \in X_\scn} \nX{x - z_\scn}.
	\end{align}
\end{theorem}
Employing the \textit{approximate Riesz representative} $R_\scn^{-1} f \in Y_\scn^\scDelta$, $\iY{R_\scn^{-1} f}{y^\scDelta_\scn} = \dY{f}{y^\scDelta_\scn}$ for all $y^\scDelta_\scn \in Y_\scn^\scDelta$, yields the \emph{localized and practical} a~posteriori error estimator
\begin{align}
	\label{approxAPostErrorEstim}
	\eta_\scn(\tilde{x}_\scn) \coloneqq \nYDp{B \tilde{x}_\scn - f} = \nY{T_\scn \tilde{x}_\scn - R_\scn^{-1} f},
\end{align}
which, under the assumptions of \cref{thmwellposedPracticalSchr}, fulfills for \textit{any} $\tilde{x}_\scn \in X_\scn$ (see \cite[Theorem 6.4]{DPGActaNumerica})
\begin{align}
	\label{apostApproxReliab}
	\frac{1}{\gamma} \, \eta_\scn(\tilde{x}_\scn) \leq \nX{x - \tilde{x}_\scn} &\leq \frac{1}{\beta} \bigl(\nu \, \eta_\scn(\tilde{x}_\scn) + \operatorname{osc}(f)\bigr), \\
	\label{oscf}
	\operatorname{osc}(f) \coloneqq \nYp{f \circ (I - \Pi)} &\leq \gamma \, \|I - \Pi\| \inf_{z_\scn \in X_\scn} \nX{x-z_\scn}.
\end{align}
The practical method admits the following equivalent mixed formulation \cite[Theorems 6.1, 6.3]{DPGActaNumerica}, which is the basis of our implementation.
\begin{proposition}
	\label{propMixedResidPract}
	Under the assumptions of \cref{thmwellposedPracticalSchr}, let $e^\scDelta_\scn \coloneqq R_\scn^{-1} f - T_\scn x_\scn \in Y_\scn^\scDelta$. Then, $x_\scn \in X_\scn$ solves \eqref{approxidealDPG} if and only if $(e^\scDelta_\scn,x_\scn)$ solves
	\begin{equation}
		\label{mixedFormDiscrete}
		\begin{alignedat}{2}
			\iY{e^\scDelta_\scn}{y^\scDelta_\scn} + b(x_\scn,y^\scDelta_\scn)& = f(y^\scDelta_\scn) &\quad &\forall \, y^\scDelta_\scn \in Y_\scn^\scDelta \\
			b(z_\scn,e^\scDelta_\scn)& = 0 &\quad &\forall \, z_\scn \in X_\scn,
		\end{alignedat}
	\end{equation}
	and in this case $\nY{e^\scDelta_\scn} = \nYDp{f - Bx_\scn} = \inf_{z_\scn \in X_\scn} \nYDp{f - Bz_\scn}$.
\end{proposition}
The mixed system thus provides the error representation function $e^\scDelta_\scn$, and with it the estimator $\eta_\scn$, as a \emph{built-in} by-product. Its algebraic realization is discussed in \cref{rmkStaticCond}.

\section*{Acknowledgments}
We thank Jay Gopalakrishnan for the helpful exchange and for providing foundational code for the primal DPG method applied to the stationary Schrödinger equation.

\section*{Declarations}
The authors used the AI assistant Claude (Anthropic) to edit and condense the text and to assist with revising and documenting the code for the numerical experiments and figures. All mathematical results and their proofs, the numerical experiments and the interpretation of their results were developed and verified by the authors. The authors assume responsibility for all content.

\bibliographystyle{siamplain}
\bibliography{references}
\end{document}